\documentclass[graybox]{svmult}

\usepackage{type1cm}        
\usepackage{makeidx}         
\usepackage{graphicx}        
\usepackage{multicol}        
\usepackage[bottom]{footmisc}
\usepackage[numbers]{natbib}
\usepackage{orcidlink}

\usepackage{newtxtext}       %
\usepackage[varvw]{newtxmath}       

\usepackage{amsmath,graphicx,textcomp,hyperref,tabularx,xcolor,theorem}

\hypersetup{
	colorlinks=true, 
	linkcolor=blue, 
	urlcolor=red, 
	citecolor=magenta,
	linktoc=all 
}

\usepackage{algorithmicx}
\usepackage{algorithm}
\usepackage{algpseudocode}

\usepackage{xpatch}
\makeatletter
\xpatchcmd{\algorithmic}
{\ALG@tlm\z@}{\leftmargin\z@\ALG@tlm\z@}
{}{}
\makeatother

\newcommand{\mathd}{\mathrm{d}}

\newenvironment{tmparmod}[3]{\begin{list}{}{\setlength{\topsep}{0pt}\setlength{\leftmargin}{#1}\setlength{\rightmargin}{#2}\setlength{\parindent}{#3}\setlength{\listparindent}{\parindent}\setlength{\itemindent}{\parindent}\setlength{\parsep}{\parskip}} \item[]}{\end{list}}

\newcounter{tmcounter}

\usepackage[T3,T1]{fontenc}
\DeclareSymbolFont{tipa}{T3}{cmr}{m}{n}
\DeclareMathAccent{\invbreve}{\mathalpha}{tipa}{16}

\newcommand{\myvector}[1]{\text{{\hspace{0.1em}}#1{\hspace{0.15em}}}}

\newcommand{\vD}{\myvector{D}}

\newcommand{\dv}[2]{{\frac{{\mathd}#1}{{\mathd}#2}}}

\newcommand{\bg}{\ensuremath{\boldsymbol{g}}}

\newcommand{\bzero}{\boldsymbol{0}}

\newcommand{\emh}{{e - \frac{1}{2}}}

\newcommand{\eph}{{e + \frac{1}{2}}}

\newcommand{\uu}{\boldsymbol{u}}
\newcommand{\uud}{\uu^{\delta}}

\newcommand{\uudhat}{\hat{\uu}^{\delta}}
\newcommand{\gnum}{\pg^{\text{avg}}}

\newcommand{\discf}{\ensuremath{\pf_h^\delta}}
\newcommand{\discfref}{\ensuremath{\hat{\pf}_{h,e}^\delta}}
\newcommand{\discfrefe}{\ensuremath{\hat{\pf}_h^\delta}}
\newcommand{\discF}{\ensuremath{\F_h^\delta}}
\newcommand{\discFref}{\ensuremath{\hat{\F}_{h,e}^\delta}}

\newcommand{\dfrx}{\partial_x^\text{FR}}
\newcommand{\dlocx}{\partial_x^\text{loc}}
\newcommand{\pdx}{\partial_x}
\newcommand{\fnum}{\pf^\text{num}}
\newcommand{\fnump}{\pf^{\text{num}+}}
\newcommand{\fnumm}{\pf^{\text{num}-}}

\newcommand{\fnumpm}{\pf^{\text{num}\pm}}

\newcommand{\fnumpj}{\pf^{(j)\text{num}+}}
\newcommand{\fnummj}{\pf^{(j)\text{num}-}}
\newcommand{\fnumpmj}{\pf^{(j)\text{num}\pm}}

\newcommand{\bB}{\boldsymbol{B}}

\newcommand{\fnumncp}{(\bB \bg)^{\text{num}^+_{\text{nc}}}}
\newcommand{\fnumncm}{(\bB \bg)^{\text{num}^-_{\text{nc}}}}
\newcommand{\fnumncpm}{(\bB \bg)^{\text{num}^\pm_{\text{nc}}}}

\newcommand{\uU}{\boldsymbol{U}}

\newcommand{\pf}{\boldsymbol{f}}

\newcommand{\fp}{\pf^+}
\newcommand{\fm}{\pf^-}

\newcommand{\pg}{\ensuremath{\boldsymbol{g}}}

\newcommand{\F}{\boldsymbol{F}}

\newcommand{\Ftotp}{\F^{\text{tot}+}}
\newcommand{\Ftotm}{\F^{\text{tot}-}}
\newcommand{\Ftotpm}{\F^{\text{tot}\pm}}

\newcommand{\ftotpj}{\pf^{(j)\text{tot}+}}
\newcommand{\ftotmj}{\pf^{(j)\text{tot}-}}
\newcommand{\ftotpmj}{\pf^{(j)\text{tot}\pm}}

\newcommand{\ftotp}{\pf^{\text{tot}+}}
\newcommand{\ftotm}{\pf^{\text{tot}-}}
\newcommand{\ftotpm}{\pf^{\text{tot}\pm}}
\newcommand{\Bloc}{\mathcal{B}^\delta}
\newcommand{\paragraphtoc}[1]{}
\providecommand{\citep}{}
\providecommand{\citet}{}

\newcommand{\poly}{\mathbb{P}}

\newcommand{\half}{\frac{1}{2}}

\newcommand{\imh}{{i - {\half}}}
\newcommand{\iph}{{i + {\half}}}

\newcommand{\Fnum}{\F^\text{num}}
\newcommand{\Fnump}{\F^{\text{num}+}}
\newcommand{\Fnumm}{\F^{\text{num}-}}
\newcommand{\Fnumpm}{\F^{\text{num}\pm}}

\newcommand{\Fnumncpm}{\F^{\text{num}^\pm_{\text{nc}}}}

\newcommand{\Fncpm}{\F^{{\text{nc}}\pm}}

\newcommand{\nc}{M}

\usepackage{url}

\makeindex             

\begin{document}

\title*{Compact Runge-Kutta flux reconstruction methods with entropy- and/or kinetic energy-preserving fluxes}
\titlerunning{Compact Runge-Kutta flux reconstruction with EC/KEP fluxes}
\author{
Arpit Babbar\orcidID{0000-0002-9453-370X},\\
Qifan Chen\orcidID{0009-0003-0790-6972}, and \\
Hendrik Ranocha\orcidID{0000-0002-3456-2277}
}
\authorrunning{Babbar, Chen, Ranocha}
\institute{
Arpit Babbar \at Institute of Mathematics, Johannes Gutenberg University, Mainz, \email{ababbar@uni-mainz.de}
\and Qifan Chen \at Department of Mathematics, The Ohio State University, Columbus, \email{chen.1234@osu.edu}
\and Hendrik Ranocha \at Institute of Mathematics, Johannes Gutenberg University, Mainz, \email{hendrik.ranocha@uni-mainz.de}
}
\maketitle
\abstract*{
Compact Runge-Kutta (cRK) methods are a class of high order methods for solving hyperbolic conservation laws characterized by their compact stencil including only immediate neighboring finite elements.
A Compact Runge-Kutta flux reconstruction (cRKFR) method for solver hyperbolic conservation laws was introduced in [Babbar, A., Chen, Q., Journal of Scientific Computing, 2025] which uses a time average flux formulation to perform evolution using a single numerical flux computation at each step, making it a single stage method.
EC/KEP numerical fluxes are often used for construction of high-order entropy stable or KEP methods for hyperbolic conservation laws, and are known to enhance the robustness of numerical methods for under-resolved simulations.
In this work, we show how these fluxes can be incorporated into the cRKFR framework for general hyperbolic equations that consist of fluxes and non-conservative products.
We test the effectiveness of this new class of methods through numerical experiments for the compressible Euler equations, magnetohydrodynamics (MHD) equations and multi-ion MHD equations.
It is observed that the application of EC/KEP fluxes enhances the robustness of the cRKFR methods.
}

\section{Introduction}
Single-stage methods solve time-dependent problems by performing evolution from one time level to the next using a single inter-element communication step, in contrast to multi-stage methods such as Runge-Kutta (RK) methods that require multiple inter-element communication steps.
These methods are of interest for modern hardware architectures where communication is often the bottleneck for performance.
Single-stage methods are an active area of research for hyperbolic equations, with the Lax-Wendroff methods~\cite{Qiu2003,Qiu2005,babbar2022lax} and ADER schemes~\cite{dumbser2008,dumbser2008a,Gaburro2023} being the most widely studied classes.
Compact Runge-Kutta (cRK) methods in a discontinuous Galerkin (DG) framework were introduced in~\cite{chen2024} as RK methods whose stencil includes only immediate neighboring elements.
These methods were extended to the flux reconstruction (FR) framework~\cite{Huynh2007} in~\cite{babbar2025crk} and to hyperbolic equations with non-conservative products in~\cite{babbar2025crknoncons}.
These compact Runge-Kutta flux reconstruction (cRKFR) methods use a time-average flux formulation to perform evolution using a single numerical flux computation at each step, making them single-stage methods like the Lax-Wendroff and ADER methods.
Their similarity to RK methods gives them advantages like the ability to treat stiff source terms using IMplicit-EXplicit (IMEX) Runge-Kutta methods, which is currently not possible for Lax-Wendroff and ADER methods.
In this work, we exploit their similarity to RK methods further by showing how entropy-conservative (EC) or kinetic energy-preserving (KEP) fluxes~\cite{tadmor1987numerical,tadmor2003entropy,jameson2008formulation,kennedygruber2008,ismail2009,chandrashekar2013,ranocha2017} can be incorporated into the cRKFR framework.
This is the first time EC/KEP fluxes have been incorporated into a single-stage method for hyperbolic equations.
This is important because these fluxes are a key ingredient in the construction of high-order entropy-stable or KEP methods, which are known to enhance robustness for under-resolved simulations \cite{chan2022,Gassner2016,sjogreen2018high,rojas2021robustness}.
The extension of cRKFR to EC/KEP fluxes is demonstrated for general hyperbolic systems with non-conservative products.
We test the effectiveness of this new class of methods through numerical experiments for the compressible Euler equations, magnetohydrodynamics (MHD) equations, and multi-ion MHD equations.

The rest of the paper is organized as follows.
In Section~\ref{sec:crkfr}, we review the cRKFR scheme for hyperbolic equations with non-conservative products as introduced in~\cite{babbar2025crknoncons}.
In Section~\ref{sec:entropy.cons.fluxes}, we show how EC/KEP fluxes can be incorporated into the cRKFR framework.
Section~\ref{sec:numerical.results} presents numerical experiments to demonstrate the improved robustness.
Finally, we conclude in Section~\ref{sec:conclusions} with a summary and future work.

\section{Review of compact Runge-Kutta flux reconstruction} \label{sec:crkfr}
Next, we review the cRKFR scheme for hyperbolic equations with non-conservative products of the form introduced in~\cite{babbar2025crknoncons}, which is an extension of the original cRKFR for conservation laws introduced in~\cite{babbar2025crk}.
We consider systems of the form
\begin{equation}
\uu_t + \pf (\uu)_x + \bB (\uu)  \bg(\uu)_x = \bzero,
\label{eq:flux.and.non.cons}
\end{equation}
where $\pf'(\uu) + \bB(\uu) \bg'(\uu)$ is a diagonalizable matrix with real eigenvalues.
The scheme was described assuming $\bg(\uu) = \uu$ in~\cite{babbar2025crknoncons}, but we present it here for the more general case of $\bg(\uu)$ as the EC fluxes are constructed for this form.
The scheme of~\cite{babbar2025crknoncons} can also handle stiff source terms, but we omit them here for simplicity.
The finite element notations, which follow the earlier works~\cite{babbar2022lax,babbar2025crk}, are introduced next.
We denote $\Omega$ as the physical domain, and $\{\Omega_e \}$ as the set of finite elements of $\Omega$ given by $\Omega_e = [x_{\emh}, x_\eph]$ so that the grid spacing is given by $\Delta x_e = x_\eph - x_{\emh}$.
The reference element is the interval $[0, 1]$ and the mapping from the physical element $\Omega_e$ to the reference element is given by $x \mapsto \xi = (x - x_{\emh})/\Delta x_e$.
The numerical solution to~\eqref{eq:flux.and.non.cons} is approximated in the basis $V_h = \{v_h : v_h |_{\Omega_e} \in \poly_N \}$, where $\poly_N$ is the space of degree $N \geq 0$, so that the solution is approximated by piecewise polynomials of degree $N$ on each element, and the solution is allowed to be discontinuous across element interfaces.
The representative of $v_h \in V_h$ is given by $\hat{v}_{h, e} = \hat{v}_{h, e} (\xi)$ defined as $\hat{v}_{h, e} (\xi) = v_h  (x_{\emh} + \xi \Delta x_e) = v_h (x)$, where $e$ is the element index which is often suppressed for brevity.
The construction of a degree $N$ basis requires $N + 1$ distinct \textit{solution points} denoted as $0 \le \xi_0 < \xi_1 < \cdots < \xi_N \le 1$, which will be Gauss-Lobatto-Legendre (GLL) nodes in this work.
Following the notations from our previous works~\cite{babbar2022lax,babbar2025crk}, the nodes and weights are with respect to the interval $[0, 1]$ whereas they are usually defined for the interval $[-1,1]$.
Denoting $\nc$ as the number of conservative variables~\eqref{eq:flux.and.non.cons}, the numerical solution $\uud_h \in V_h^{\nc}$, inside an element $\Omega_e$ is given in reference coordinates as $\uudhat_{h, e} (\xi, t) = \sum_{p = 0}^N \uu_{e, p} (t) \ell_p (\xi)$, where each $\ell_p$ is a Lagrange polynomial of degree $N$ defined to satisfy $\ell_p (\xi_q) = \delta_{p q}$.
The spatial derivatives required to solve~\eqref{eq:flux.and.non.cons} are computed on the reference interval $[0, 1]$ using the differentiation matrix $\vD = [D_{pq}]$ whose entries are given by
\begin{equation}
D_{pq} = \ell_q' (\xi_p), \qquad 0 \le p, q \le N.
\label{eq:diff.matrix}
\end{equation}
\subsection{FR discretization for the conservative flux}
A crucial ingredient of the cRKFR scheme and the FR scheme in general is the degree $N$ \textit{discontinuous flux} approximation $\discf$
, which is defined in reference coordinates for element $e$ as
\begin{equation} \label{eq:discts.flux}
\discfref(\xi) = \sum_{p=0}^N \pf(\uu_{e,p}) \ell_p(\xi).
\end{equation}
Then, the FR differentiation operator is defined as
\begin{equation}
\begin{gathered}
\dfrx \pf (\uud_h) = \pdx \pf_h, \quad \pf_h = \discf +
(\fnum_\eph - \fm_\eph) g_R + (\fnum_{\emh} - \fp_\emh) g_L, \\
\fm_\eph = \discfrefe (1), \quad \fp_\emh = \discfrefe (0),
\label{eq:dfrx.defn}
\end{gathered}
\end{equation}
where $g_L, g_R \in \poly_{N + 1}$ are FR correction functions~\cite{Huynh2007,Vincent2011a} satisfying $g_L (0) = g_R (1) = 1$, $g_L (1) = g_R (0) = 0$, and $\fnum_\eph$ is the numerical flux at the interface $x_\eph$. We use the Rusanov flux~\cite{rusanov1962}
\begin{equation} \label{eq:num.flux.fr}
\fnum_\eph(\uud_h) = \half(\pf(\uu_\eph^-) + \pf(\uu_\eph^+)) - \frac {\lambda_\eph}{2} (\uu_\eph^+ - \uu_\eph^-), \quad \uu_\eph^\pm = \uud_h(x_\eph^\pm),
\end{equation}
where $\lambda_\eph = \lambda(\uu_\eph^-,\uu_\eph^+)$ is a Rusanov/local Lax-Friedrichs~\cite{rusanov1962} estimate of the wave speed at the interface $x_\eph$, which for the case $\bB = \bzero$, is given by
\begin{equation} \label{eq:rusanov.wave.speed}
\lambda(\uu_l, \uu_r) = \max_{\uu \in \{\uu_l, \uu_r\}} \sigma(\pf'(\uu)),
\end{equation}
where $\sigma(A)$ denotes the spectrum of the matrix $A$.
Since $\pf_h$ defined in~\eqref{eq:dfrx.defn} is globally continuous (taking the numerical flux value $\fnum_{e + 1 / 2}$ at interface $e + 1 / 2$), it is called the \textit{continuous
flux approximation} in the FR literature~{\cite{Huynh2007}}.
Following~\cite{chen2024,babbar2025crk}, we define a local differentiation operator $\dlocx$ as
\begin{equation}
\dlocx \pf (\uud_h) (\xi_p) = \pdx \discf (\xi_p) .
\label{eq:dlocx.defn}
\end{equation}

\subsection{Finite volume discretization for non-conservative product}
Next, we explain the discretization of the non-conservative term by describing the finite volume method (FVM) for equation~\eqref{eq:flux.and.non.cons} without the conservative flux, i.e.,
\begin{equation}
\uu_t = - \bB (\uu)  \pg(\uu)_x. \label{eq:pure.non.conservative}
\end{equation}
As discussed in~\cite{babbar2025crknoncons}, a natural discretization of~\eqref{eq:pure.non.conservative} is obtained by using a central difference approximation
\begin{equation}
\dv{\uu_i}{t} = - \frac{1}{\Delta x}  (\fnumncm_\iph - \fnumncp_\imh), \label{eq:semi.fv}
\end{equation}
where $\fnumncpm_\iph = \fnumncpm(\uu_i, \uu_{i+1})$ which is defined as
\begin{equation}
\fnumncpm(\uu_-, \uu_+) = \half \bB({\uu_\pm}) \gnum(\uu_-, \uu_+) = \half \bB({\uu_\pm}) (\pg(\uu_-) + \pg(\uu_+)).
\label{eq:non.cons.num.flux.basic}
\end{equation}
The crucial difference between~\eqref{eq:non.cons.num.flux.basic} and the standard numerical fluxes for conservation laws is that each finite volume interface $\iph$ now has two numerical fluxes $\fnumncpm_\iph$, which is equivalent to saying that the numerical fluxes are no longer continuous across the interfaces.
From the common experience with conservation laws, the above scheme will need to have additional dissipation terms to be stable.
Thus, the final numerical flux for solving~\eqref{eq:pure.non.conservative} will actually need to be chosen as
\begin{equation}
\fnumncpm_\iph
- \frac{\lambda_\iph}{2}  (\uu_{i + 1} - \uu_i),
\label{eq:num.flux.non.cons}
\end{equation}
where $\lambda_\iph = \lambda(\uu_i, \uu_{i+1})$ is defined by the Rusanov wave speed estimate~\eqref{eq:rusanov.wave.speed} extended to non-conservative products as $\lambda(\uu_l, \uu_r) = \max_{\uu \in \{\uu_l, \uu_r\}} \sigma(\bB(\uu) \pg'(\uu))$.
However, unlike the central part~\eqref{eq:non.cons.num.flux.basic}, the dissipation term is continuous across the interface $\iph$.
Hence, we include the dissipation term of the entire system in the conservative numerical flux.
Thus, for the system with both conservative and non-conservative terms~\eqref{eq:flux.and.non.cons}, we define the numerical fluxes $\fnumpm_\iph$ as
\begin{equation}
\begin{split}
\fnumpm_\iph &=
\fnumpm(\uu_i, \uu_{i+1}) :=
\fnum(\uu_i, \uu_{i+1}) + \fnumncpm(\uu_i, \uu_{i+1}), \\
\fnum_\iph &= \fnum(\uu_i, \uu_{i+1}) :=
\half(\pf(\uu_i) + \pf(\uu_{i+1}))
- \frac {\lambda_\iph}{2} (\uu_{i+1} - \uu_{i}),
\end{split}
\label{eq:combined.num.flux}
\end{equation}
where the dissipation coefficient is defined as
\begin{equation}
\lambda_\iph = \max_{\uu \in \{\uu_i, \uu_{i+1}\}} \sigma(\pf'(\uu) + \bB(\uu) \pg'(\uu)).
\label{eq:combined.wave.speed}
\end{equation}

\subsection{cRKFR scheme for the complete non-conservative equation}
We now describe the cRKFR scheme for solving~\eqref{eq:flux.and.non.cons} using the above numerical fluxes.
The general FR~\cite{Huynh2007} scheme to solve~\eqref{eq:flux.and.non.cons} that applies to any choice of correction functions and solution points is derived via DG methods in Appendix~B of~\cite{babbar2025crknoncons}.
In this work, we always use GLL solution points and $g_2$ correction functions, because this choice gives semi-discrete entropy stability for the standard RKFR scheme when used with EC fluxes~\cite{fisher2013,carpenter2014entropy,Gassner2016,gassner2016shallow,ranocha2017}.
However, we still give the general description for the most part in this section.
Here we directly state the semi-discretization of~\cite{babbar2025crknoncons}:
\begin{multline}
\dv{\uud_h}{t}
+ \dlocx \pf(\uud_h) + \bB (\uud_h) (\dlocx \bg(\uud_h)) \\
+ [
  \fnumm_\eph - \ftotm_\eph
  ] \partial_x g_R  + [
  \fnump_{\emh} - \ftotp_{\emh}
  ] \partial_x g_L = \bzero
,
\label{eq:general.non.conservative.fr.one}
\end{multline}
where $\dlocx$ is the local derivative operator~\eqref{eq:dlocx.defn}
, $\fnumpm_\eph$ is defined using the finite volume conservative and non-conservative numerical fluxes~(\ref{eq:combined.num.flux}, \ref{eq:non.cons.num.flux.basic}) as
\begin{equation}
\begin{split}
\fnumpm_\eph &= \fnum_\eph + \fnumncpm_\eph, \qquad \fnumncpm_\eph = \fnumncpm(\uu_\eph^-, \uu_\eph^+),\\
\fnum_\eph &= \fnum(\uu_\eph^-, \uu_\eph^+),
\end{split} \label{eq:combined.num.flux.fr}
\end{equation}
where $\uu_\eph^{\pm}$ are the extrapolations as in~\eqref{eq:num.flux.fr}.
The $\ftotpm_\eph$ consist of the extrapolations of the conservative and non-conservative terms at the element interface $x_\eph$ given by
\begin{equation} \label{eq:total.flux.non.cons}
\ftotpm_\eph = \bB (\uu_\eph^\pm) \bg(\uu_\eph)^\pm + \pf_\eph^\pm,
\end{equation}
where, as in~\eqref{eq:dfrx.defn}, $\pf_\eph^\pm$ are extrapolations obtained from the discontinuous fluxes~\eqref{eq:discts.flux} at the interface $x_\eph$.
The motivation for why the extrapolations of the non-conservative terms are chosen as in~\eqref{eq:total.flux.non.cons} is given in Appendix~B of~\cite{babbar2025crknoncons}.
It is further shown there that the scheme~\eqref{eq:general.non.conservative.fr.one} reduces to the first-order finite volume scheme~\eqref{eq:combined.num.flux} when $N = 0$, and that for the FR scheme with GLL solution points and $g_2$ correction function can be written in terms of the FR operators $\dfrx$ like in the conservative case as
\begin{equation}
\begin{gathered}
\dv{\uud_h}{t} +  \dfrx \pf(\uud_h) + \bB (\uud_h) \dfrx  \bg(\uud_h) = \bzero,\qquad
\dfrx \bg(\uud_h) := \partial_x  \bg(\uu)_h,
\\ \bg(\uu)_h := \bg(\uud_h) +
(\bg^\text{num}_\eph - \bg(\uu_\eph)^-)  g_R +
(\bg^\text{num}_{\emh} - \bg(\uu_{\emh})^+)  g_L.
\end{gathered} \label{eq:fr.non.cons.gll}
\end{equation}
In order to now derive the cRKFR scheme for~\eqref{eq:flux.and.non.cons}, we show the fully discrete scheme obtained by applying an explicit RK method to~\eqref{eq:general.non.conservative.fr.one}.
An $s$ stage explicit RK method is specified by its Butcher tableau
\begin{equation}\label{eq:butcher}
\begin{array}{c|c}
{c}&{A}\\\hline
&{b}\\
\end{array},\quad {A} = ({a}_{ij})_{s\times s}, \quad {b} = ({b}_1,\dots, {b}_s), \quad {c}_i = \sum_j {a}_{ij},
\end{equation}
where ${A}$ is a lower triangular matrix in \eqref{eq:butcher}.
We obtain the fully discrete RKFR scheme for the non-conservative system~\eqref{eq:flux.and.non.cons} to be
\begin{subequations}\label{eq:rkfr.non.cons}
\begin{align}
&\begin{aligned}
\ \ \uu^{(i)}& = \uu^n -  \Delta t\sum_{j = 1}^{i-1} {a}_{ij} \Big( \dlocx \pf(\uu^{(j)}) + \bB (\uu^{(j)}) \dlocx \bg(\uu^{(j)}) \\
&\qquad\qquad\qquad + \big[\fnummj_\eph - \ftotmj_\eph\big]\partial_x g_R
+ \big[\fnumpj_{\emh} - \ftotpj_{\emh}\big]\partial_x g_L \Big),
\\
& \qquad\qquad\qquad \qquad i = 1, 2, \dots, s,
\end{aligned}\label{eq:rkfr1} \\
&\begin{aligned}
\uu^{n+1} &= \uu^n - \Delta t \sum_{j = 1}^s {b}_j \Big( \dlocx \pf(\uu^{(j)} ) + \bB (\uu^{(j)}) \dlocx \bg(\uu^{(j)}) \\
&\qquad\qquad\qquad + \big[\fnummj_\eph - \ftotmj_\eph\big]\partial_x g_R
+ \big[\fnumpj_{\emh} - \ftotpj_{\emh}\big]\partial_x g_L \Big),
\end{aligned}
\label{eq:rkfr2}
\end{align}
\end{subequations}
where $\fnumpmj_\eph = \fnumpm(\uu_\eph^{(j)-}, \uu_\eph^{(j)+})$ are the sum of conservative and non-conservative numerical fluxes~\eqref{eq:combined.num.flux.fr} at stage $j$ evaluated as~\eqref{eq:combined.num.flux}, and $\ftotpmj_\eph$ are the extrapolations of the conservative and non-conservative terms at stage $j$ computed as in~\eqref{eq:total.flux.non.cons}.
The cRKFR scheme in the time-average framework~\cite{babbar2025crk} is obtained by dropping the inter-element terms from the inner stages~\eqref{eq:rkfr1} while using FR operators on the time averages in the final update~\eqref{eq:rkfr2} as follows
\begin{eqnarray}
\uu^{(i)} & = & \uu^n -  \Delta t\sum_{j = 1}^{i-1} {a}_{ij} \Big( \dlocx \pf(\uu^{(j)}) + \bB (\uu^{(j)}) \dlocx \bg(\uu^{(j)}) \Big),
\;\; i = 1, \ldots, s, \label{eq:crkfr.inner}\\
\uu^{n + 1} & = & \uu^n
- \Delta t \pdx \discF
- \Delta t \Bloc_h \nonumber \\
&& \quad - \Delta t (\Fnumm_\eph - \Ftotm_\eph)\pdx g_R - \Delta t (\Fnump_\emh - \Ftotp_\emh) \pdx g_L,
\label{eq:crkfr.non.cons.final}
\\ \Ftotpm_\eph &=& \F_\eph^\pm  + \Fncpm_\eph,\quad \Fnumpm_\eph = \Fnum_\eph + \Fnumncpm_\eph, \label{eq:total.F.non.cons}
\end{eqnarray}
where the local time averaged flux $\discF$ in~\eqref{eq:crkfr.non.cons.final} is computed following~\cite{babbar2025crk,babbar2025crknoncons} as
\begin{equation}
  \discFref = \sum_{p = 0}^N \F_{e, p} \ell_p (\xi), \qquad \F_{e, p} =
  \sum_{i = 1}^s {b}_i  \pf (\uu_{e, p}^{(i)}), \label{eq:disc.avg.flux}
\end{equation}
the local time average of the non-conservative product, $\Bloc_h$, is computed as a time average of the local approximations of $\bB(\uu) \partial_x \bg(\uu)$ used in~(\ref{eq:fr.non.cons.gll},~\ref{eq:rkfr1},~\ref{eq:rkfr2}) to give
\begin{equation}
\Bloc_h = \sum_{j=1}^s {b}_j \bB (\uu^{(j)}) \dlocx \bg(\uu^{(j)}).
\label{eq:local.time.avg.non.cons}
\end{equation}
The conservative part of the time averaged numerical flux $\Fnum_\eph$ is~\cite{babbar2025crk}
\begin{equation} \label{eq:numflux.cons}
\Fnum_\eph = \half (\F_\eph^- + \F_\eph^+) - \frac{\lambda_\eph}{2} (\uU_\eph^+ - \uU_\eph^-),
\end{equation}
with the dissipation coefficient $\lambda_\eph$ using the wave speed estimate for the full non-conservative system~\eqref{eq:flux.and.non.cons} as in~\eqref{eq:combined.wave.speed}.
Using the time-average solution $\uU_{\eph}^\pm$ (computed similar to~\eqref{eq:disc.avg.flux}) in the dissipative term instead of the extrapolated solution $\uu_\eph^\pm$ enhances the stability of the scheme~\cite{babbar2022lax,babbar2025crk}.
The conservative part $\F_\eph^\pm$ in the extrapolation term $\Ftotpm_\eph$ in~(\ref{eq:total.F.non.cons},~\ref{eq:numflux.cons}) is computed by extrapolating the time averaged flux~\eqref{eq:disc.avg.flux} to the element interface $x_\eph$ similar to~\eqref{eq:total.flux.non.cons}.
The non-conservative part $\Fncpm_\eph$ of the extrapolation term~$\Ftotpm_\eph$ in~\eqref{eq:total.F.non.cons} is computed as $\Fncpm_\eph = \sum_j {b}_j \bB (\uu_\eph^{(j)\pm}) \bg(\uu_\eph^{(j)\pm}).$
The expression of the non-conservative time averaged numerical flux $\Fnumncpm_\eph$ that is obtained by time averaging~(\ref{eq:combined.num.flux.fr},~\ref{eq:non.cons.num.flux.basic}) directly would be
\begin{equation}
\Fnumncpm_\eph = \sum_{j=1}^s {b}_j \fnumncpm_\eph = \sum_{j=1}^s {b}_j \bB (\uu_\eph^{(j), \pm})  \bg^\text{avg}_\eph(\uu^{(j)}).
\label{eq:bad.numflux.non.cons}
\end{equation}
However, as pointed out in~\cite{babbar2025crknoncons}, equation ~\eqref{eq:bad.numflux.non.cons} requires the inter-element communication of the inner stage values $\{\uu_\eph^{(j), \pm} \}_{j = 1}^s$ at the element interfaces between neighboring elements.
This would be similar to the communication volume required for the first cRKDG scheme proposed in~\cite{chen2024}, and the ADER scheme~\cite{dumbser2008} as ADER requires communication of the local space-time predictor at the element interfaces when dealing with non-conservative products.
Thus, the increased communication volume just might be unavoidable for some practical problems.
However, in the large class of problems tried in~\cite{babbar2025crknoncons} and the additional problems tried in this work, inter-element communication has been reduced by using $\Fnumncpm_\eph = \bB (\uu_\eph^{n\pm})  \bg^\text{avg}_\eph(\uU)$.
This computation does not require any additional communication compared with the conservative case as the communication of $\uU_\eph^\pm$ is already needed in~\eqref{eq:numflux.cons}.
It has been used in all the numerical results of this work and in~\cite{babbar2025crknoncons}, and has given optimal order of accuracy for problems with smooth analytical solutions and results that are extremely close to those obtained by~\eqref{eq:bad.numflux.non.cons} for all tested problems.

\section{cRKFR with entropy- or kinetic energy-preserving fluxes} \label{sec:entropy.cons.fluxes}

We propose to use EC/KEP fluxes in the cRKFR scheme in a fashion similar to the standard DG/FR schemes as described, for example, in~\cite{Gassner2016,gassner2016shallow}.
This is often referred to as the \emph{flux differencing form} of the DG/FR schemes~\cite{lefloch2002fully,fisher2013,carpenter2014entropy,Gassner2016}.
In this form, the application of the derivative operator $\dlocx$ defined in~\eqref{eq:dlocx.defn} applied to the point-wise discontinuous flux~\eqref{eq:discts.flux} and to the non-conservative term in~\eqref{eq:general.non.conservative.fr.one} is replaced by standard EC/KEP fluxes for the conservative and non-conservative terms.
The element interface fluxes $\fnumpm_\eph$ are still computed using the Rusanov flux~\eqref{eq:combined.num.flux.fr} and thus the overall scheme still has a dissipative behavior.
The first step in our description is to write the application of the differentiation matrix~\eqref{eq:diff.matrix} in the ``inner stages'' of the cRKFR scheme~\eqref{eq:crkfr.inner} explicitly as
\begin{align}
(\dlocx \pf(\uud_h))(\xi_p) &= \frac{1}{\Delta x} \sum_{q=0}^N D_{pq} \pf_q, \ \ \pf_q = \pf(\uud_h(\xi_q)), \label{eq:flux.dev}\\
(\bB (\uud_h) \dlocx \bg(\uud_h))(\xi_p) &= \frac{1}{\Delta x} \sum_{q=0}^N D_{pq} \bB_p \bg_q, \ \  \bB_p = \bB (\uud_h(\xi_p)), \ \bg_q = \bg(\uud_h(\xi_q)) \label{eq:noncons.dev}.
\end{align}
The flux differencing form for the conservative case~\eqref{eq:flux.dev} is given by~\cite{fisher2013}
\begin{equation}
(\dlocx \pf(\uud_h))(\xi_p) = \frac{1}{\Delta x} \sum_{q=0}^N 2D_{pq} \fnum_\text{S}(\uu_p, \uu_q),
\label{eq:entropy.cons.flux.dev}
\end{equation}
where $\fnum_\text{S}$ is a symmetric flux satisfying $\fnum_\text{S}(\uu_p, \uu_p) = \pf(\uu_p)$.
Equation~\eqref{eq:entropy.cons.flux.dev} can be seen as a generalization of~\eqref{eq:flux.dev} by taking $\fnum_\text{S}(\uu_p, \uu_q) = (\pf_p + \pf_q)/2$.
For the compressible Euler equations, several choices of the flux function $\fnum_\text{S}$ are available to obtain entropy conservation~\cite{ismail2009,chandrashekar2013, ranocha2017}.
In addition to being EC, some of these fluxes~\cite{chandrashekar2013,ranocha2017} also preserve kinetic energy.
However, these fluxes are often expensive to compute as they involve logarithmic mean evaluations.
Thus, a cheaper alternative to add robustness is to use fluxes that only preserve kinetic energy~\cite{kennedygruber2008} by using only arithmetic averages.
Similarly, the non-conservative fluxes are introduced by replacing~\eqref{eq:noncons.dev} by
\begin{equation}
(\bB (\uud_h) \dlocx \bg(\uud_h))(\xi_p) = \frac{1}{\Delta x} \sum_{q=0}^N D_{pq} (\bB \bg)_{\text{NS}}^\text{num}(\uu_p, \uu_q),
\label{eq:entropy.cons.noncons.dev}
\end{equation}
where $(\bB \bg)_{\text{NS}}^\text{num}$ is a non-symmetric function.
Equation~\eqref{eq:entropy.cons.noncons.dev} generalizes~\eqref{eq:noncons.dev} as we can choose $(\bB \bg)_{\text{NS}}^\text{num}(\uu_p, \uu_q) = \bB (\uu_p) \bg(\uu_q)$.
EC fluxes have been constructed for several equations including the shear shallow water equations~\cite{yadav2023}, ideal magnetohydrodynamics (MHD)~\cite{ramirez2024,bohm2020} and for multi-ion MHD~\cite{ramirez2025}.
In~\cite{Artiano2026,Waruszewski2022}, more specific forms of $(\bB \bg)_{\text{NS}}^\text{num}$ are analyzed and general procedures to construct such fluxes to obtain entropy conservation is given.

The formulas~(\ref{eq:flux.dev},~\ref{eq:noncons.dev}) directly apply to the inner stages of the cRKFR scheme~\eqref{eq:crkfr.inner}.
In addition, in the final update~\eqref{eq:crkfr.non.cons.final}, the computation of the local flux derivative $\pdx \discF$ and the non-conservative term $\Bloc_h$~\eqref{eq:local.time.avg.non.cons} are also modified to be
\begin{align*}
(\pdx \discF)_p & = \frac{1}{\Delta x} \sum_{j=1}^s {b}_j \sum_{q=0}^N 2 D_{pq} \fnum_\text{S}(\uu_p^{(j)}, \uu_q^{(j)}),\\
(\Bloc_h)_p &= \frac{1}{\Delta x} \sum_{j=1}^s {b}_j \sum_{q=0}^N D_{pq} (\bB \bg)_{\text{NS}}^\text{num}(\uu_p^{(j)}, \uu_q^{(j)}).
\end{align*}
For the standard FR scheme, such a modification of the flux derivative and the non-conservative term is sufficient to obtain semi-discrete entropy or kinetic energy preservation.
For the cRKFR scheme, there is no semi-discretization and thus such a statement cannot be made.
However, the numerical results show an improvement of robustness when using the EC fluxes in the cRKFR scheme, as is observed in the standard FR schemes (see e.g.,~\cite{chan2022}).

\section{Numerical results} \label{sec:numerical.results}

As described in Section~\ref{sec:entropy.cons.fluxes}, we use EC/KEP fluxes for the volume terms in the cRKFR scheme and the dissipative Rusanov flux~\eqref{eq:num.flux.fr} at the interfaces.
Using EC/KEP fluxes is known to improve the robustness of numerical schemes with RK time integration.
In this section, we demonstrate this improvement in robustness also for the cRKFR scheme for high-gradient problems for the compressible Euler equations, ideal MHD equations and multi-ion MHD equations~\cite{chan2022,ramirez2025}.
It is consistently observed that the EC/KEP fluxes allow the cRKFR scheme to be stable for significantly longer times than the standard cRKFR scheme.
We use the EC flux from~\cite{ranocha2017} and the KEP flux from~\cite{kennedygruber2008} for the Euler equations, and the fluxes from~\cite{ramirez2024,bohm2020} and~\cite{ramirez2025} for the ideal MHD equations and multi-ion MHD equations, respectively.
We only show the results for compressible Euler equations and multi-ion MHD equations here; similar improvements are observed for the ideal MHD equations.
The physical fluxes for the compressible Euler equations can be found in our previous works~\cite{babbar2024admissibility,babbar2025crk}, and the fluxes and non-conservative terms for the multi-ion MHD equations can be found in~\cite{ramirez2025}.
They are not shown here to save space.
The code to reproduce the results of this work is available at~\cite{babbar2026crkes}, including all the initial conditions and equations.
We use the Julia \cite{bezanson2017julia} package \texttt{Tenkai.jl}~\cite{tenkai}, which further uses \texttt{Trixi.jl}~\cite{schlottkelakemper2021purely,schlottkelakemper2020trixi,ranocha2022adaptive} through Julia's package manager to access the EC flux functions.

\subsection{Convergence test}

\begin{figure}[htb]
\centering
\begin{tabular}{cc}
\includegraphics[width=0.44\columnwidth]{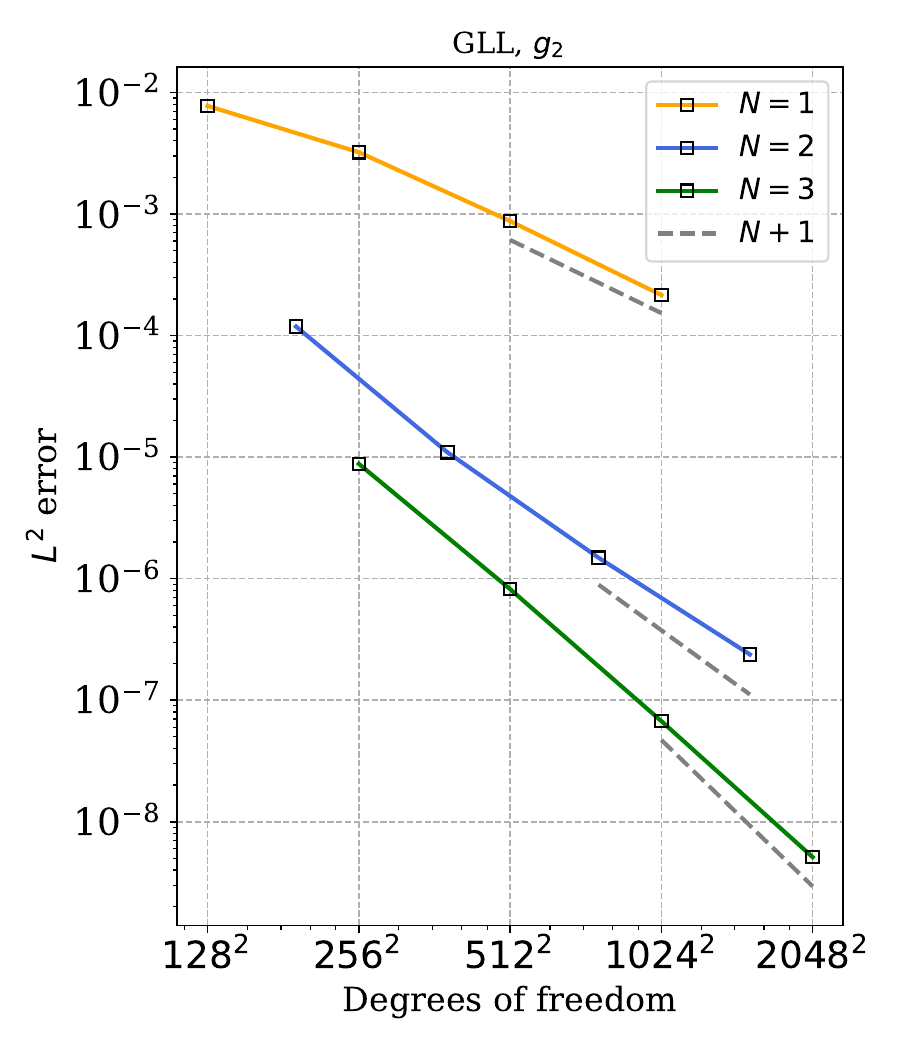} &
\includegraphics[width=0.44\columnwidth]{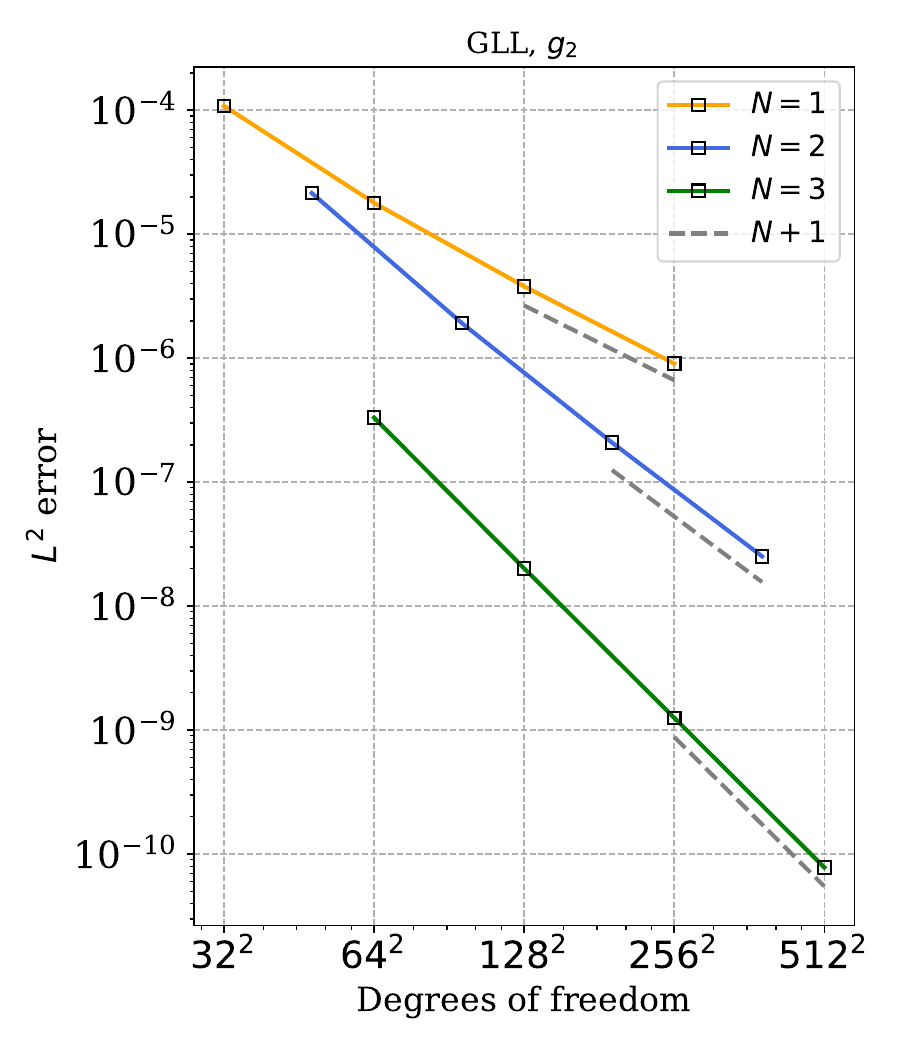} \\
(a) & (b)
\end{tabular}
\caption{Convergence test for the cRKFR scheme for (a) the isentropic vortex for the compressible Euler equations (with KEP flux from~\cite{kennedygruber2008}) and (b) manufactured solution for the multi-ion MHD equations (with the EC flux from~\cite{ramirez2025}).}
\label{fig:convergence}
\end{figure}

For the first convergence test, we use the cRKFR scheme with EC~\cite{ranocha2017} and KEP~\cite{kennedygruber2008} fluxes for the isentropic vortex test~\cite{Yee1999,Spiegel2016} for the compressible Euler equations.
The initial condition can be found in~\cite{babbar2024admissibility,babbar2025crk} or in our reproducibility repository~\cite{babbar2026crkes} and is not shown here to save space.
Since both the EC flux and the KEP flux nearly give the same accuracy, we only show the results with the KEP flux~\cite{kennedygruber2008} in Figure~\ref{fig:convergence}a.
Optimal order of accuracy is not observed for this problem, which is a known issue when using GLL nodes~\cite{babbar2022lax}.
The accuracy is same the as with standard cRKFR scheme of~\cite{babbar2025crk}.
The second convergence test is for the multi-ion MHD equations and is constructed by the method of manufactured solutions as described in~\cite{ramirez2025}.
The initial conditions and source terms for this problem can be found in~\cite{ramirez2025} or in our reproducibility repository~\cite{babbar2026crkes} and are not shown here to save space.
Figure~\ref{fig:convergence}b shows results using the EC fluxes from~\cite{ramirez2025}, demonstrating optimal order of accuracy.

\subsection{Kelvin-Helmholtz instabilities}

\begin{figure}[htb]
\centering
\includegraphics[width=0.6\columnwidth]{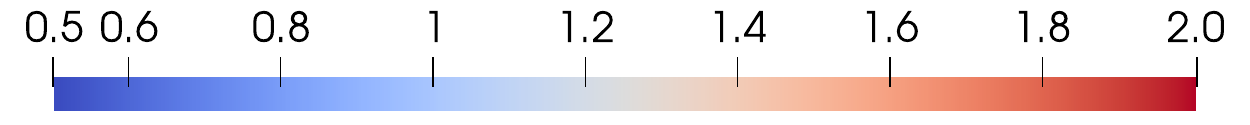} \\
\begin{tabular}{ccc}
\includegraphics[width=0.32\columnwidth]{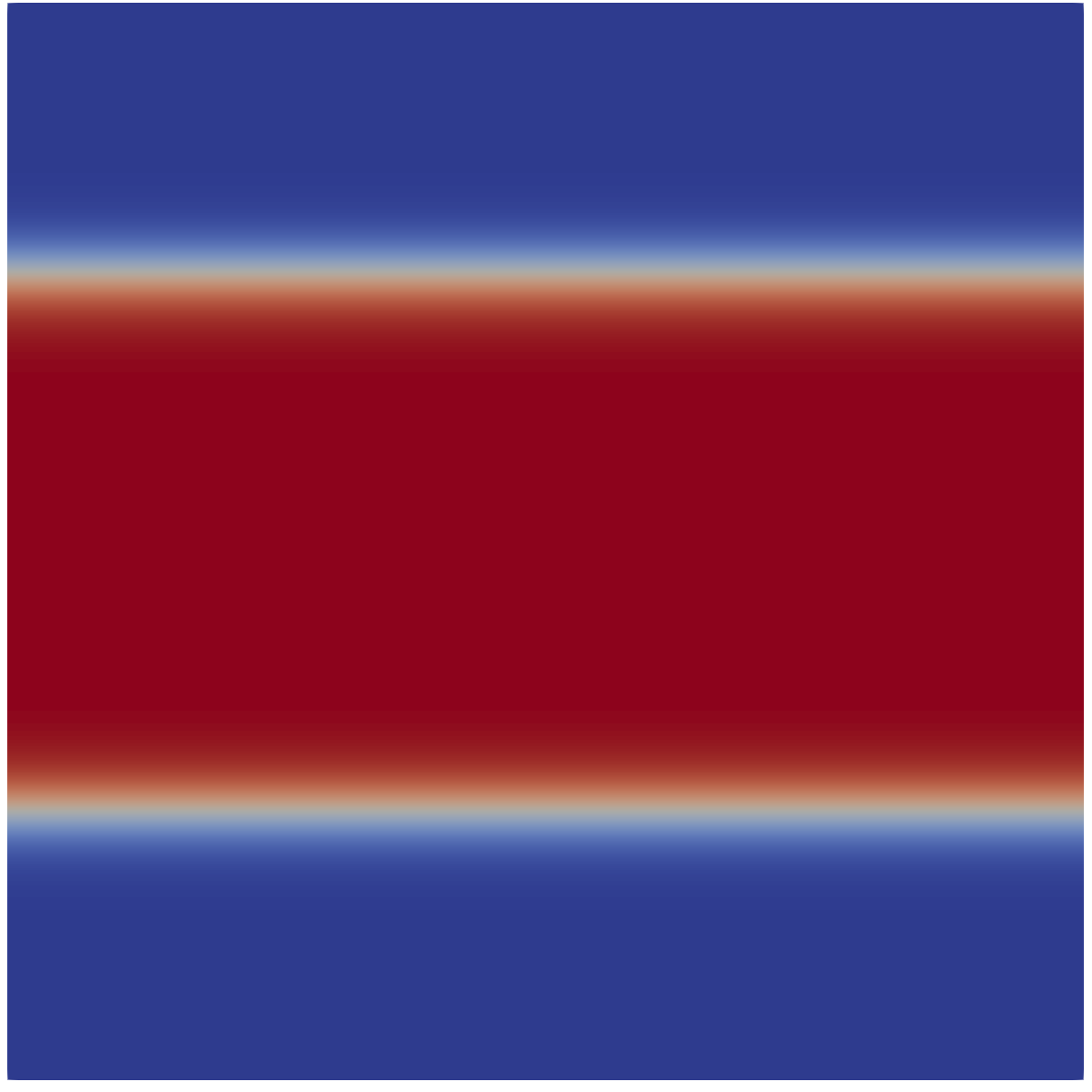} &
\includegraphics[width=0.32\columnwidth]{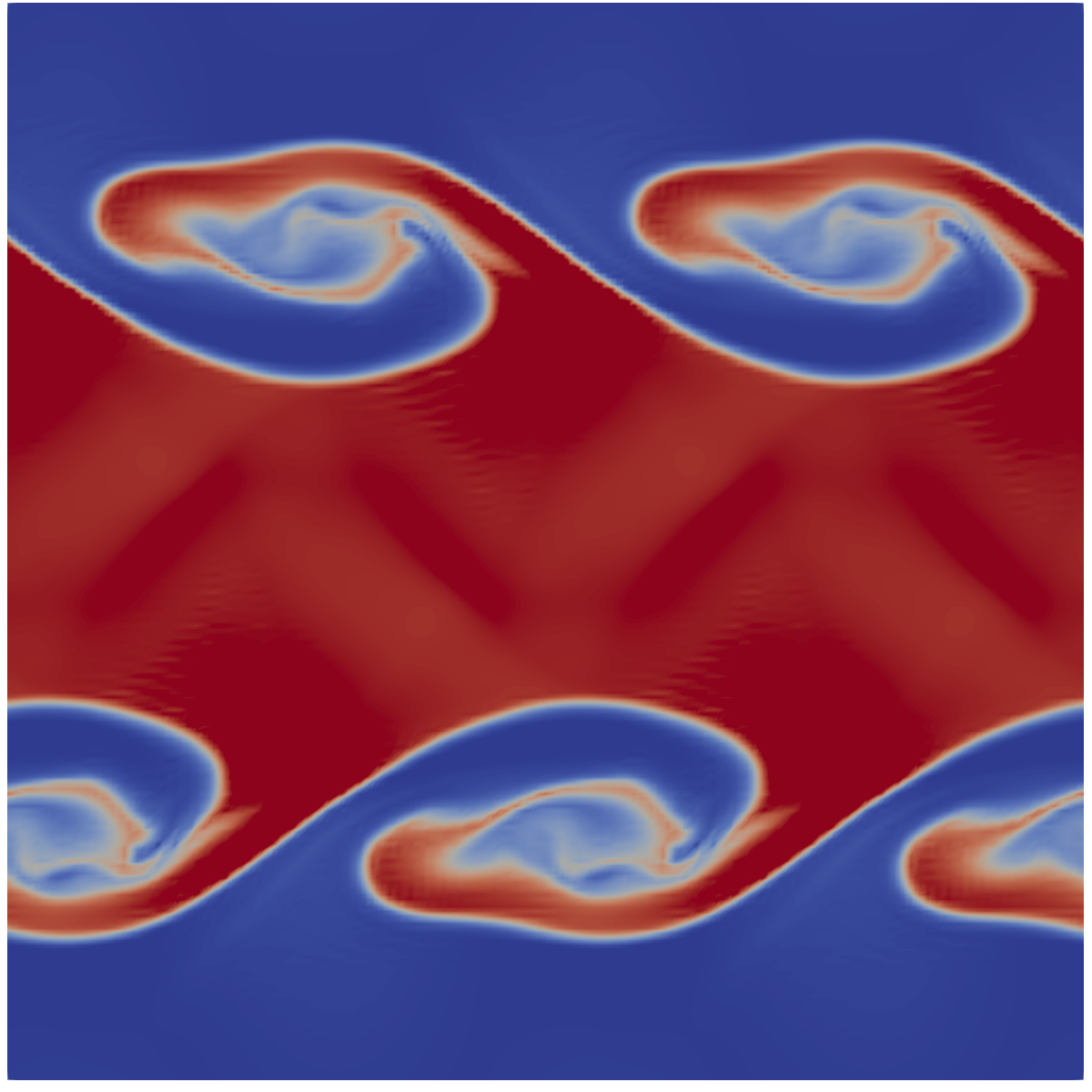} &
\includegraphics[width=0.32\columnwidth]{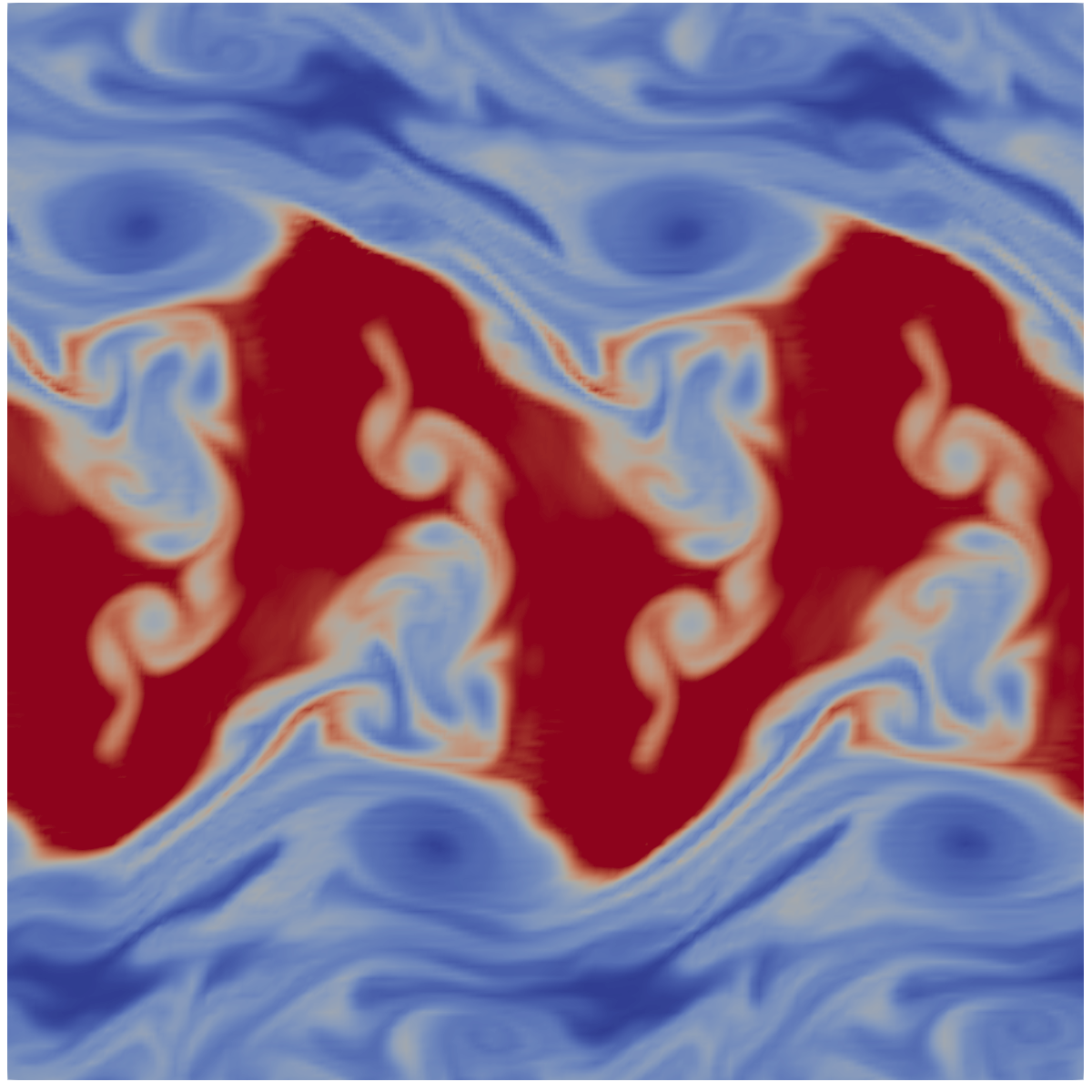} \\
$t=0$ & $t=3$ & $t=15$
\end{tabular}
\caption{Kelvin-Helmholtz instability for the compressible Euler equations using the cRKFR scheme with the KEP flux of~\cite{kennedygruber2008} and the blending scheme of~\cite{babbar2024admissibility,babbar2025crk} with $\alpha_\text{max} = 0.008$.
The density profile is shown at different times using polynomial degree $N=3$ and $64 \times 64$ elements.}
\label{fig:khi}
\end{figure}

\begin{figure}[htb]
\centering
\includegraphics[width=0.6\columnwidth]{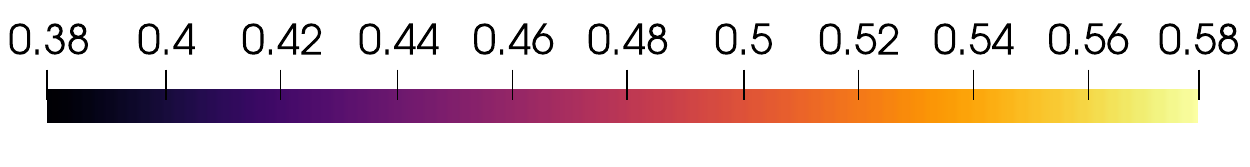} \\
\begin{tabular}{cccc}
\includegraphics[width=0.21\columnwidth]{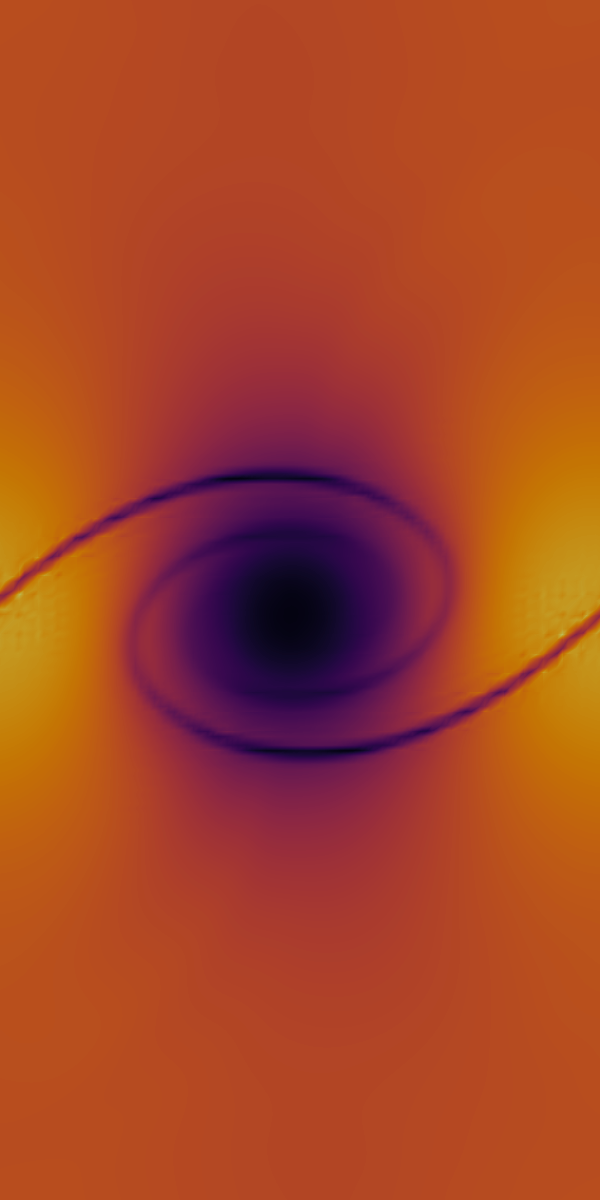} &
\includegraphics[width=0.21\columnwidth]{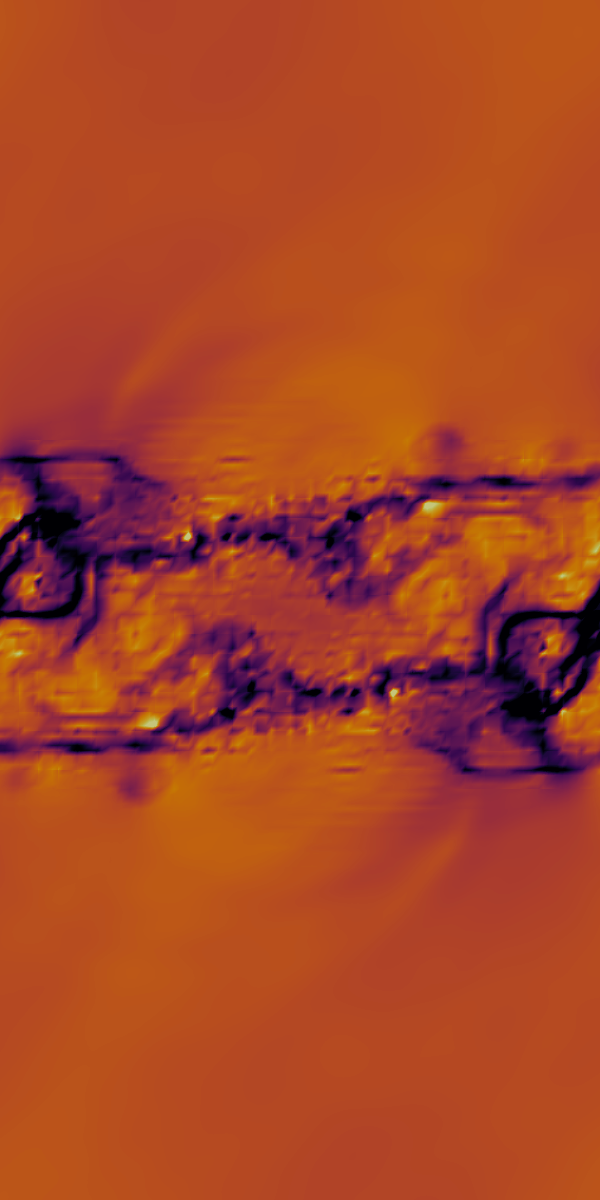} &
\includegraphics[width=0.21\columnwidth]{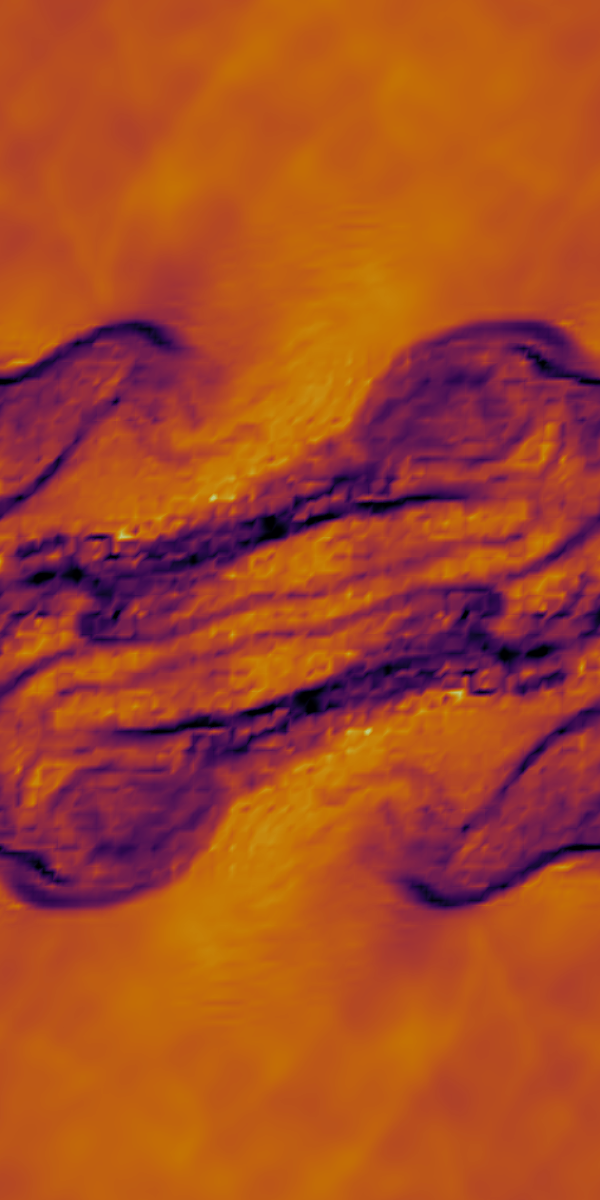} &
\includegraphics[width=0.21\columnwidth]{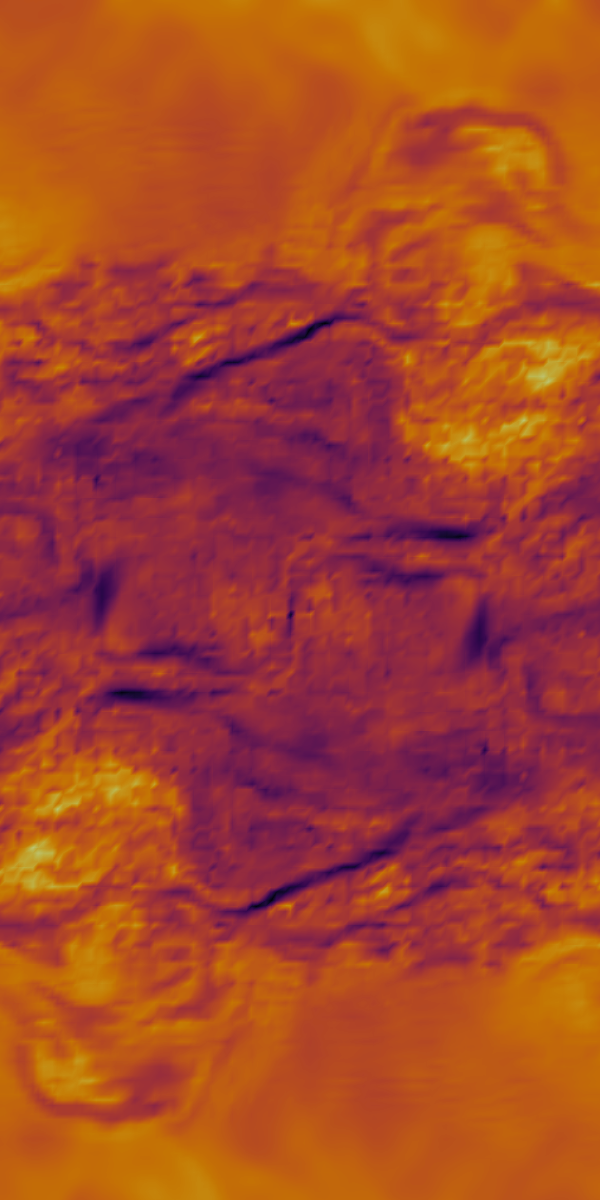} \\
$t=5$ & $t=8$ & $t=12$ & $t=20$
\end{tabular}
\caption{Multi-ion MHD Kelvin-Helmholtz instability using the EC flux of~\cite{ramirez2025} without the subcell-based blending limiter.
The density profile of ion species 1 is shown at different times using polynomial degree $N=3$ and $64 \times 128$ elements.}
\label{fig:khi.multi.ion}
\end{figure}

Kelvin-Helmholtz instabilities occur in various physical applications such as fluid mechanics, atmospheric physics and astrophysics, and are often used to test robustness.
We test the scheme for the Kelvin-Helmholtz instability for the compressible Euler equations and multi-ion MHD equations.
In the case of the compressible Euler equations, the initial condition is taken from~\cite{ramirez2021proceedings} and is explicitly given by $(\rho, v_1, v_2, p) = (1/4 + 3 B /4, (B-1)/2, \sin (2 \pi x)/10, 1)$ where $B(x,y) = \tanh(15y + 7.5) - \tanh(15y-7.5)$ on the domain $[-1, 1] \times [-1, 1]$ with periodic boundary conditions.
The scheme required the blending limiter of~\cite{babbar2024admissibility,babbar2025crk} to be stable for this problem.
With the EC flux of~\cite{ranocha2017} or the KEP flux of~\cite{kennedygruber2008}, a blending coefficient limited by $\alpha_\text{max} = 0.008$ was sufficient to stabilize the scheme for this problem till $t=15$.
The fact that we are able to obtain the same robustness with the inexpensive KEP flux of~\cite{kennedygruber2008} is a positive result for the practical usage.
However, without the EC/KEP fluxes, the scheme crashes around $t=4$ using the same $\alpha_\text{max}$ value.
The results with the KEP flux are shown in Figure~\ref{fig:khi} for various time levels.
For the multi-ion MHD equations, the Kelvin-Helmholtz instability setup is taken from~\cite{ramirez2025} and is not shown here to save space.
The scheme with the EC flux of~\cite{ramirez2025} was observed to be stable without any limiter for this problem, and the results are shown in Figure~\ref{fig:khi.multi.ion} for various time levels till $t=20$.
Without the EC flux, the scheme crashes around $t=5$ for this problem.

\subsection{Richtmyer-Meshkov instability}

\begin{figure}[htb]
\centering
\includegraphics[width=0.8\columnwidth]{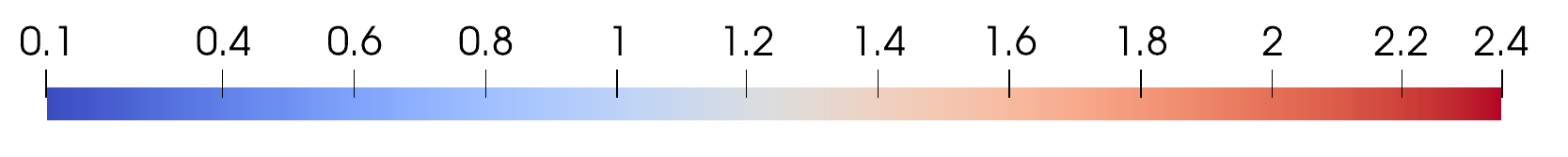}\\
\begin{tabular}{cccc}
\includegraphics[width=0.21\columnwidth]{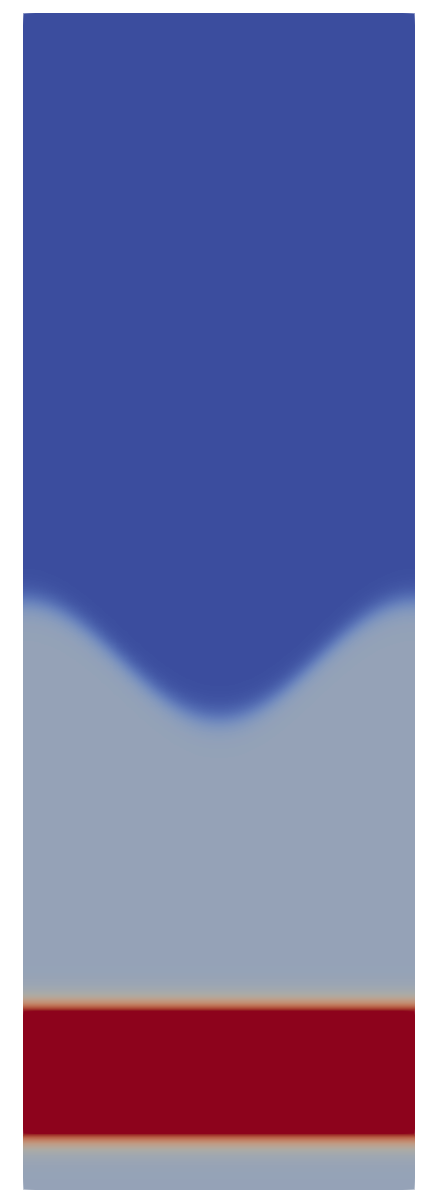} &
\includegraphics[width=0.21\columnwidth]{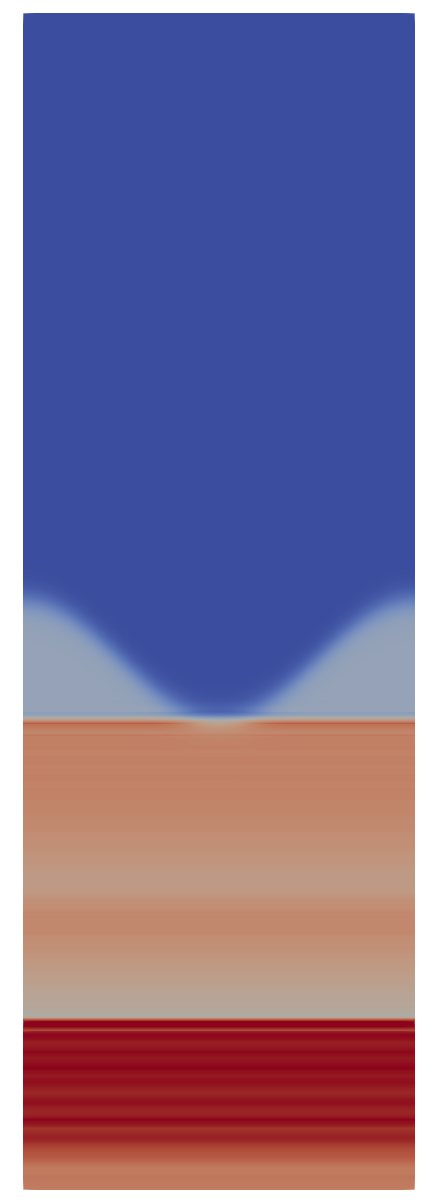} &
\includegraphics[width=0.21\columnwidth]{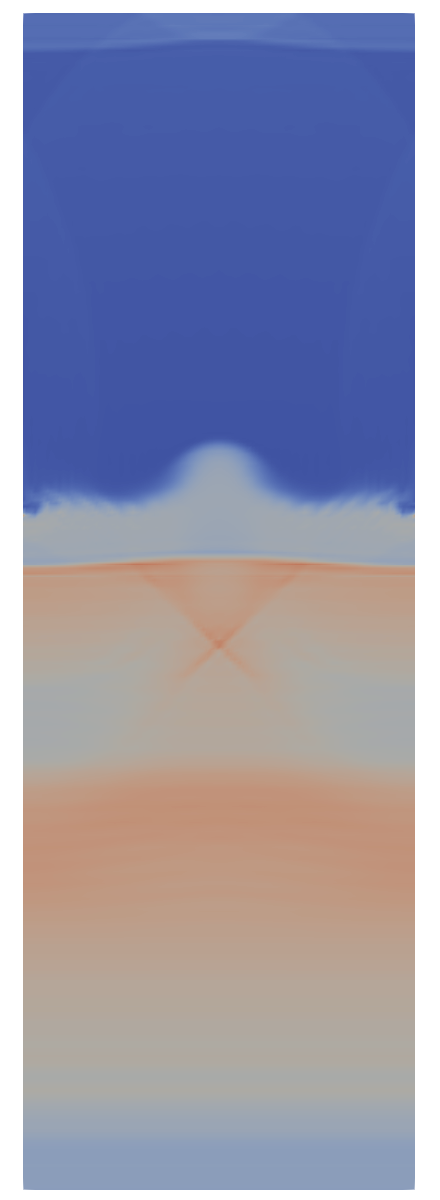} &
\includegraphics[width=0.21\columnwidth]{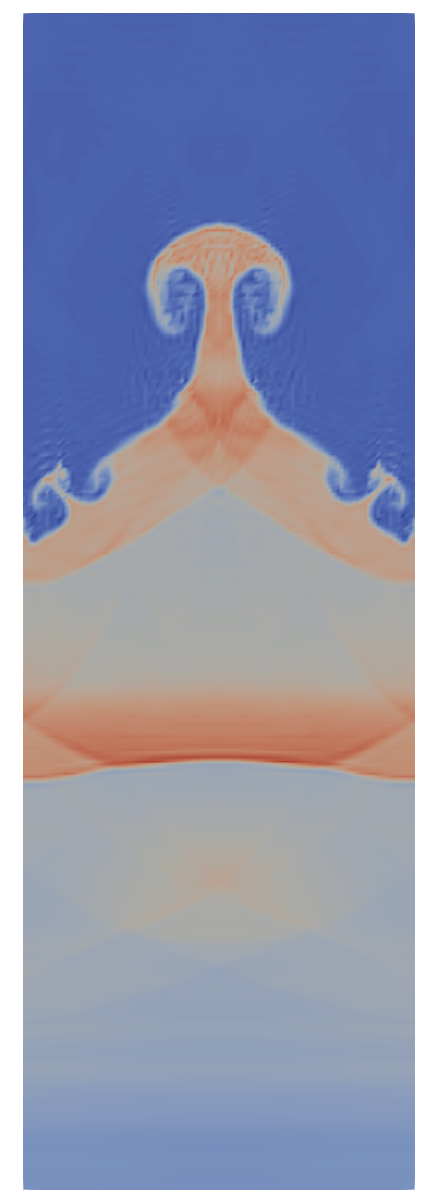} \\
$t=0$ & $t=6$ & $t=15$ & $t=30$
\end{tabular}
\caption{Richtmyer-Meshkov instability for the compressible Euler equations using the cRKFR scheme with the KEP flux of~\cite{kennedygruber2008} and blending scheme of~\cite{babbar2024admissibility,babbar2025crk} with $\alpha_\text{max} = 0.001$. The density profile is shown at different times using polynomial degree $N=3$ and $32 \times 96$ elements.}
\label{fig:richtmyer.meshkov}
\end{figure}

The Richtmyer-Meshkov instability generates small-scale flow features by passing a shock over a stratified fluid~\cite{Richtmyer1960,Meshkov1972}.
This test on the physical domain $[0, 40/3] \times [0,40]$ is also taken from~\cite{chan2022}, but we provide the initial condition as there is a typo in~\cite{chan2022}.
Let $d_{a,b} = a + (1 + \tanh(sx))(b-a)/2$ denote a smooth approximation to a discontinuous function that transitions from $a$ to $b$ at $x=0$ with slope $s = 2$.
The initial density is $\rho(x,y) = d_{1, 1/4}(y - (18 + 2 \cos (6 \pi x/L))) + d_{3.22,0}(|y-4| - 2)$ and the initial pressure is $p(x,y) = d_{4.9, 1}(|y-4| - 2)$ with $L = 40$.
The scheme required the blending limiter of~\cite{babbar2024admissibility,babbar2025crk} to be stable for this problem.
The cRKFR scheme with EC flux of~\cite{ranocha2017} or the KEP flux of~\cite{kennedygruber2008} and blending coefficient limited by $\alpha_\text{max} = 0.001$ was stable, and could be run until $t=30$.
However, without the EC/KEP fluxes, the scheme crashes around $t=8.7$ with the same $\alpha_\text{max}$ value.
The results with the KEP fluxes are shown in Figure~\ref{fig:richtmyer.meshkov} for various time levels.

\section{Conclusions} \label{sec:conclusions}

EC and KEP fluxes have been incorporated into the volume terms of the time average compact Runge-Kutta flux reconstruction (cRKFR) scheme of~\cite{babbar2025crk,babbar2025crknoncons}.
The time average framework makes the cRKFR scheme a single-stage method requiring only one inter-element numerical flux evaluation per time step.
The proposed method is applicable to conservative systems and systems with non-conservative products, with EC fluxes used for both terms.
Numerical experiments show that the EC and KEP fluxes improve the robustness of the cRKFR scheme.
The experiments were performed for the compressible Euler equations and multi-ion MHD equations, and optimal order of accuracy is obtained for problems with smooth analytical solutions.
This is the first time that EC fluxes have been used with a single-stage high-order method.
This work gives numerical evidence that EC fluxes can be used in a single-stage method to improve the robustness of the scheme.
An important area for future research will be provable entropy stability for single-stage methods, e.g., by using relaxation Runge-Kutta methods.

\section*{Acknowledgments}

AB and HR were supported by the Deutsche Forschungsgemeinschaft
(DFG, German Research Foundation, project number 528753982
as well as within the DFG priority program SPP~2410 with project number 526031774).
AB was also supported by the Alexander von Humboldt Foundation.

\bibliographystyle{spbasic}
\bibliography{references.bib}
\end{document}